\documentclass[12pt,leqno]{article}

\usepackage{amssymb,amsmath,amsthm}
\usepackage{mathrsfs}
\numberwithin{equation}{section}

\newcommand{\corps}{\operatorname{\mathbf{k}}}

\newcommand{\Int}{{\rm Int}}

\newcommand{\supp}{{\rm supp}}

\newcommand{\lan}{\langle}
\newcommand{\ran}{\rangle}
\newcommand{\eps}{\varepsilon}
\newcommand{\Supp}{\operatorname{Supp}}
\newcommand{\dist}{{\rm dist}}
\newcommand{\Ss}{\operatorname{SS}}
\newcommand{\codim}{\operatorname{codim}}
\providecommand{\bysame}{\rule{2.5em}{.1mm}\,}
\newcommand{\indlim}[1][]{\mathop{\varinjlim}\limits_{#1}}

\def\etens{\mathbin{\boxtimes}}
\newcount\temp
\temp=\catcode`\@
\catcode`\@=11
\def\@bletens{\mathbin{\etens^{L}}}
\def\@letens_#1{\mathbin{\etens_{\raise1.5ex\hbox to-.1em{}#1}^{L}}}
\def\letens{\@ifnextchar _{\@letens}{\@bletens}}
\catcode`\@=\temp

\def\phi{{\varphi}}
\def\epsilon{\varepsilon}

\newcommand{\ba}{\begin{array}}
\newcommand{\ea}{\end{array}}

\def\shb{\mathscr{B}}

\def\shd{\mathscr{D}}
\def\she{\mathcal{E}}

\def\shm{\mathscr{M}}
\def\shn{\mathscr{N}}
\def\sho{\mathscr{O}}

\newcommand{\C}{\mathbb{C}}

\newcommand{\R}{\mathbb{R}}

\newcommand{\Z}{\mathbb{Z}}

\newcommand{\isoto}[1][]{\xrightarrow[#1]{\sim}}

\newcommand{\tens}[1][]{\otimes_{#1}}

\newcommand{\oim}[1]{{#1}_*}
\newcommand{\opb}[1]{#1^{-1}}
\renewcommand{\hom}[1][]{{\mathcal{H}\kern-1pt{om}_{\raise1.5ex\hbox to.1em{}#1}}}
\newcommand{\roim}[1]{{R#1}_*}
\newcommand{\reim}[1]{{R#1}_!}

\newcommand{\rhom}[1][]{{R{\mathcal{H}}{om}_{\raise1.5ex\hbox to.1em{}#1}}}
\newcommand{\ext}[1][]{{\mathcal{E}xt}_{\raise1.5ex\hbox to.1em{}#1}}

\newcommand{\rsect}{R\Gamma}

\newcommand{\Hom}[1][]{\mathrm{Hom}_{\raise1.5ex\hbox to.1em{}#1}}
\newcommand{\RHom}[1][]{\mathrm{RHom}_{\raise1.5ex\hbox to.1em{}#1}}

\def\car{{\operatorname{Ch}}}

\theoremstyle{plain}

\newtheorem{theorem}{Theorem}[section]
\newtheorem{proposition}[theorem]{Proposition}
\newtheorem{lemma}[theorem]{Lemma}
\newtheorem{corollary}[theorem]{Corollary}

\theoremstyle{definition}

\newtheorem{definition}[theorem]{Definition}

\newtheorem{example}[theorem]{Example}
\newtheorem{examples}[theorem]{Examples}
\newtheorem{remark}[theorem]{Remark}

\newcommand{\eq}{\begin{eqnarray}}
\newcommand{\eneq}{\end{eqnarray}}
\newcommand{\eqn}{\begin{eqnarray*}}
\newcommand{\eneqn}{\end{eqnarray*}}

\newenvironment{nnum}{
  \begin{enumerate}
  \itemsep=0pt
  
  }
  {\end{enumerate}}

\newenvironment{anum}{
  \begin{enumerate}
  \itemsep=0pt
  
  }
  {\end{enumerate}}

\newcommand{\bnum}{\begin{nnum}}
\newcommand{\enum}{\end{nnum}}
\newcommand{\banum}{\begin{anum}}
\newcommand{\eanum}{\end{anum}}
\newcommand{\ol}{\overline}
\def\dddt{{\raise-.3em\hbox{$\big\cdot$}}}

\newcommand{\cl}{\colon}

\newcommand{\hot}{\hookleftarrow}
\newcommand{\reg}{{\operatorname{reg}}}
\newcommand{\bl}{\bigl}
\newcommand{\br}{\bigr}
\newcommand{\Ker}{{\operatorname{Ker}}}
\newcommand{\Coker}{{\operatorname{Coker}}}

\begin{document}

\author{Masaki Kashiwara,
Teresa Monteiro Fernandes%
\\ 
and Pierre Schapira}
\title{Truncated microsupport and holomorphic solutions of
 D-Modules}
\date{}

\maketitle
\footnote{The research of the second author
was supported by FCT and Programa Ci{\^e}ncia,
Tecnologia e Inova\c c{\~a}o do Quadro Comunit{\'a}rio de Apoio.}

\begin{abstract}
We study the truncated microsupport $\Ss_k$ of sheaves on a real manifold. 
Applying our results
to the case of $F=\rhom[\shd](\shm,\sho)$, the complex of
holomorphic solutions of a coherent $\shd$-module $\shm$,
we show that $\Ss_k(F)$ is completely
determined by
the characteristic variety of $\shm$.
As an application,
we obtain an extension theorem for the sections of $H^j(F)$, $j<d$,
defined on an open subset whose boundary is
non characteristic outside of a complex analytic 
subvariety of codimension $d$.
We also give a characterization of the perversity for
$\C$-constructible sheaves in terms of their truncated microsupports.
\end{abstract}

\section{Introduction}\label{section:intro}
The notion of microsupport of sheaves 
was introduced in the course of the study 
of the theory of linear partial differential equations (LPDE),
and it is now applied in various domains of mathematics.
References are made to \cite {K-S1}.

For an  object  $F$ of the derived category of abelian 
sheaves on a real
manifold $X$, its microsupport $\Ss(F)$ is a closed conic subset of the
cotangent bundle $\pi\cl T^*X\to X$ which describes 
the direction of ``non-propagation'' of $F$. In particular, for a
smooth closed submanifold  $Y$ of $X$, 
the support $\Supp(\mu_Y(F))$ 
of the Sato microlocalization $\mu_Y(F)$ of $F$ along $Y$ is contained in
$\Ss(F)\cap T^*_YX$, where $T^*_YX$ denotes the conormal bundle to $Y$ in
$X$.

If $X$ is a complex manifold and $\shm$ is a system
of LPDE, that is, a coherent 
module over the sheaf
$\shd_X$ of holomorphic differential operators, the complex $F$ of
holomorphic solutions of this system is given by
$\rhom[\shd_X](\shm, \sho_X)$,
and the microsupport of $F$ is then
the characteristic variety $\car(\shm)$ of $\shm$.

However some phenomena of propagation may happen 
in specific degrees, related to
the principle of unique continuation for holomorphic functions, and
this leads to a variant of the notion of microsupport, that of
 ``truncated microsupport'', a notion introduced by the authors 
of \cite{K-S1}
but never published.

Following a suggestion of these authors,  F.~Tonin \cite{To}
was able to regain in the language of the truncated microsupport 
a result of 
Ebenfelt-Khavinson-Shapiro \cite{E-K-S1,E-K-S2}. 
In their papers, these authors obtained the 
extension of holomorphic solutions of a
differential operator when the solutions are defined 
on an open subset with smooth boundary
non-characteristic outside a smooth complex hypersurface. 
Note that the problem of extending holomorphic solutions 
across non characteristic real hypersurfaces plays a crucial role in
the theory of LPDE, and was initiated by J.~Leray \cite{Le}, 
followed by \cite{Ze}, \cite{B-S}, \cite{Ka-1} (see also \cite{Ho} and
\cite{S} for an exposition of these results).

The truncated microsupport is defined as follows.
Let $\corps$ be a field, let $X$ be a real manifold and let 
$F\in D^b(\corps_X)$ be an object of the
derived category of sheaves of $\corps$-vector spaces on $X$. 
For an integer $k\in\Z$,
a point $p\in T^*X$ does not belong to the truncated microsupport
$\Ss_k(F)$ if and only if
$F$ is microlocally at $p$ isomorphic to
an object of $D^{>k}(\corps_X)$. 
Hence $\{\Ss_k(F)\}_{k\in\Z}$ is
an increasing sequence, while 
$\Ss_k(F)$ is an empty set for $k\ll0$ and
$\Ss_k(F)$ coincides with $\Ss(F)$ for $k\gg0$.

In this paper, we give equivalent definitions of the
truncated microsupport and we study its 
behavior under exterior tensor product,
smooth inverse image and proper direct image. 
We introduce then the closed subset $\Ss_k^Y(F)$ of $\opb{\pi}(Y)$ 
which describes the support of the microlocalization of $F$ along
submanifolds of $Y$ and prove (see Theorem \ref{T:6}):
\eqn
\Ss_k (F)=\overline{\Ss_k (F)\setminus \opb{\pi}(Y)}\cup \Ss_k^Y(F).
\eneqn

We then apply this result to the complex 
$F=\rhom[\shd_X](\shm,\sho_X)$ of holomorphic solutions of a coherent
$\shd_X$-module $\shm$
on a complex manifold $X$. 

Let  $S$ be a closed complex analytic subset of
codimension greater than or equal to $d$
and let  $S'$ be a closed complex analytic subset of $S$
of codimension greater than $d$
such that $S_0:=S\setminus S'$ is a smooth
submanifold of codimension $d$. 
We first prove the estimate below (see Proposition \ref{T:7}):

\eqn\label{eq:1.1}
&&
\begin{array}{rl}
\Ss_{d-1}(F)&=\overline{\Ss_{d-1}(F)\setminus\pi^{-1}(S)},\\[5pt]
\Ss_{d}(F)&=\overline{\Ss_{d}(F)\setminus\pi^{-1}(S)}
\cup \ol{\Ss_{d}(F)\cap T^*_{S_0}X}.
\end{array}
\eneqn
In particular, if $j: \Omega\hookrightarrow X$  is the embedding 
of a pseudo-convex open subset with smooth
boundary $\partial \Omega$, and if  $\partial \Omega$ is transversal to $S$
(i.e.\ $T^*_{\partial \Omega}X\cap\ol{T^*_{S_0}X}\subset T^*_XX$)
and non characteristic for $\shm$ outside of
$S$, one has 
$$\she {xt}^j_{\shd_X}(\shm, \sho^{+}_X/
\sho_X)=0\quad \mbox{for any $j<d$,}
$$ 
where $\sho^{+}_X/\sho_X=j_{*}j^{-1}\sho_X/\sho_X$.

Next we calculate $\Ss_k(F)$ 
in terms of the characteristic variety of $\shm$
(see Theorem \ref{th:sskFcharM}).
Letting $\car(\shm)=\bigcup_{\alpha\in A}V_\alpha$
be the decomposition of $\car(\shm)$ into irreducible components,
one has
\eq
&&\Ss_k(F)=\Bigl(\mathop\cup\limits_{\codim \pi(V_\alpha)<k}V_\alpha\Bigr)
\cup
\Bigl(\mathop\cup\limits_{\codim \pi(V_\alpha)=k}T^*_{\pi(V_\alpha)}X\Bigr).
\eneq
In particular, if $F$ is a perverse sheaf
({\em i.e.,} $\shm$ is holonomic),
then letting
$\Ss(F)=\bigcup_{\alpha\in A} \Lambda_\alpha$
be the decomposition into irreducible components,
one has
$$\Ss_k(F)=\bigcup_{\codim \pi(\Lambda_\alpha)\le k}\Lambda_\alpha.$$
Conversely if $F\in D^b(\C_X)$ is $\C$-constructible
and if it satisfies
\eqn
\Ss_k(F)\cup\Ss_k(\rhom(F,\C_X))\subset
\bigcup_{\codim \pi(\Lambda_\alpha)\le k}\Lambda_\alpha
\quad\mbox{for every $k$},
\eneqn
then $F$ is a perverse sheaf.

\medskip
\noindent
{\sl Acknowledgment}\quad
The authors would like to thank A.~D'Agnolo for 
helpful discussion.

\section{Notations and review}\label{S:1}

We will mainly follow the notations in \cite{K-S1}.

Let $X$ be a real analytic manifold. We denote by
$\tau\colon TX\to X$ the tangent bundle to $X$ and
by $\pi\colon T^*X\to X$ the
 cotangent bundle. We identify $X$ with the zero section of $T^*X$ and set
$\dot{T}^*X=T^*X\setminus X$.  We denote by
$\dot{\pi}\colon
\dot{T}^*X\to X$  the restriction of $\pi$ to $\dot{T}^*X$.
For a smooth submanifold $Y$ of $X$, $T_YX$ denotes
 the normal bundle to $Y$ and $T^*_YX$ the conormal bundle.
In particular, $T^*_X X$ is identified with $X$.

For a submanifold $Y$ of $X$ and a subset  $S$ of $X$,
we denote by $C_Y(S)$ the normal cone to $S$
along $Y$,  a closed conic subset of $T_YX$.

For a morphism $f\colon X\to Y$ of real manifolds, we denote
by
$$f_{\pi}\colon X\times_YT^*Y\to T^*Y
\mbox{ and }f_d \colon X\times_YT^*Y\to T^*X$$
the associated morphisms.

For a subset $A$ of $T^*X$,
we denote by $A^a$ the image of $A$ by the
antipodal map $a\colon (x;\xi)\mapsto (x;-\xi)$.
The closure of $A$ is denoted by  $\overline{A}$.
For a cone $\gamma\subset TX$,
the polar cone $\gamma^\circ$ to $\gamma$ is the convex cone in $T^*X$
defined by
$$\gamma^\circ=\{(x;\xi)\in T^*X;
\mbox{$x\in\pi(\gamma)$ and $\langle v,\xi\rangle\geq 0$ for any
$(x;v)\in \gamma$}\}.$$

Let $\corps$ be a field.
We denote by $D(\corps_X)$
the derived category  of complexes of sheaves of $\corps$-vector spaces
on $X$,
and by $D^b(\corps_X)$ the full subcategory of
$D(\corps_X)$ consisting of complexes with bounded cohomologies.

For $k\in \Z$, we denote as usual by  $D^{\geq k}(\corps_X)$
(resp.\ $D^{\leq k}(\corps_X)$)
the full additive subcategory of $D^b(\corps_X)$ consisting of
objects $F$ satisfying $H^j(F)=0$ for any $j<k$
(resp.\ $H^j(F)=0$ for any $j>k$).
The category 
$D^{\ge k+1}(\corps_X)$ is sometimes denoted by
$D^{>k}(\corps_X)$.

We denote by $\tau^{\leq k}\cl D(\corps_X)\to D^{\leq k}(\corps_X)$
the truncation functor. Recall that for $F\in D(\corps_X)$
the morphism $\tau^{\le k}F\to F$
induces isomorphisms
$H^j(\tau^{\le k}F)\isoto H^j(F)$
for $j\le k$ and $H^j(\tau^{\le k}F)=0$
for $j>k$.

If $F$ is an object of $D^b(\corps_X)$,
$\Ss(F)$ denotes its microsupport, a closed
$\R^{+}$-conic involutive subset of $T^*X$.
For $p\in T^*X$, $D^b(\corps_X; p)$ denotes the localization of 
$D^b(\corps_X)$ by the full triangulated subcategory consisting of objects $F$
such that $p\notin \Ss(F)$.

If $Y$ is a submanifold,  $\mu_Y(F)$ denotes the Sato microlocalization of
$F$ along $Y$. Recall that $\mu_Y(F)\in D^b(\corps_{T_Y^*X})$ and
\begin{equation}\label{E:13}
H^j(\mu_Y ( F))_p\simeq\indlim[Z] H^j_{Z}(F)_{\pi(p)}
\quad\mbox{ for $p\in T_Y^*X$ and $j\in \Z$,}
\end{equation}
where   $Z$ runs through
the family of closed subsets of $X$ such that
\begin{equation}
 C_{Y}(Z)_{\pi(p)}\setminus\{0\}\subset
\{v\in (T_{Y}X)_{\pi(p)}; \langle v, p\rangle>0\}.
\end{equation}
On a complex manifold $X$, we consider the sheaf  $\sho_X$  of
holomorphic functions and the sheaf $\shd_X$
of linear holomorphic differential operators of finite order.
Concerning the theory of $\shd$-modules, references are made to
\cite{Ka-2}.

\section{Truncated microsupport}\label{S:2}
We shall give here several equivalent definitions of the truncated
microsupport.

For  a closed cone $\gamma \subset\R^n$, one sets  
\eqn
&& Z_{\gamma}:=\{(x;y)\in \R^n \times \R^n; x-y \in \gamma\}.
\eneqn
Let $q_1,\,q_2\colon \R^n\times \R^n\to \R^n$
be the first and the second projections.

One defines 
the integral transform with kernel
$\corps_{Z_{\gamma}}$,
\eqn
&&\corps_{Z_\gamma}\circ\ : D^b (\corps_{\R^n})\to D^b (\corps_{\R^n}),
\qquad
\corps_{Z_\gamma}\circ G = \reim{q_1}(\corps_{Z_\gamma}\otimes\opb{q_2} G).
\eneqn
If $G$ has compact support, one has the following
formula for the stalk of $\corps_{Z_{\gamma}}\circ G $ at $x\in \R^n$:
\eqn
&&(\corps_{Z_{\gamma}}\circ G)_x\simeq \rsect(\R^n; \corps_{x+\gamma ^a}
\otimes G).
\eneqn

Recall that a closed convex cone
$\gamma$ is called {\em proper}
if $0\in\gamma$ and $\Int(\gamma^\circ)\neq\emptyset$.

For $(x_0;\xi_0)\in \R^n\times(\R^n)^*$ and $\epsilon\in\R$ we set:
\eqn
&& H_{\epsilon}(x_0,\xi_0)=\{x\in \R^n;
\langle x-x_0,
\xi_0\rangle>-\epsilon\}.
\eneqn
If there is no risk of confusion,
we write $H_\epsilon$ instead of $H_{\epsilon}(x_0,\xi_0)$
for short.
The following result is proved in \cite{K-S1}.

\begin{lemma}\label{L:01}
Let $X$ be an open subset of $\R^{n}$ and let   $G\in D^b(\corps_X)$.
Let  $p=(x_0;\xi_0)\in T^*X$.
Then
$p\notin \Ss(G)$ if and only if there exist an open neighborhood $W$ of $x_0$,
a proper closed convex cone $\gamma$
and
$\epsilon>0$
such that $\xi_0\in \Int(\gamma^{\circ})$,
$(W+{\gamma^a})\cap \ol{H_{\epsilon}}\subset X$ and
$$  H^j(X; \corps_{(x+{\gamma^a})\cap H_{\epsilon}}\otimes
G)=0\quad\mbox{for any $j\in \Z$ and $x\in W$.}$$
\end{lemma}

We shall give a similar version of the above lemma for 
the truncated microsupport.

\begin{proposition}\label{L:1}
Let $X$ be a real analytic manifold and let $p\in T^*X$.
Let $F\in D^b(\corps_X)$ and $k\in\Z$, 
$\alpha\in\Z_{\ge1}\cup\{\infty,\omega\}$.
Then the following conditions are equivalent:
\bnum
\item[{\rm (i)}$_{k}$]  There exists $F'\in D^{>k}(\corps_X)$
and an isomorphism $F\simeq F'$ in $D^b(\corps_X;p)$.
\item[{\rm (ii)}$_{k}$] There exists $F'\in D^{>k}(\corps_X)$ and a
morphism $F'\to F$ in $D^b(\corps_X)$
which is an isomorphism in $D^b(\corps_X; p)$.
\item[{\rm (iii)}$_{k,\alpha}$] There exists an open conic neighborhood $U$
of $p$ such that for any
$x\in \pi(U)$ and for any $\R$-valued $C^{\alpha}$-function $\phi$
defined on a neighborhood of $x$ such
that $\phi(x)=0$, $d\phi(x)\in U$, one has
\begin{equation}\label{E:1}
H^j_{\{\phi\geq 0\}}(F)_x=0\quad\mbox{for any $j\le k$.}
\end{equation}
\enum
When $X$ is an open subset of $\R^n$ and $p=(x_0;\xi_0)$,
the above conditions are also equivalent to
\bnum
\item[{\rm (iv)}$_{k}$]
There exist a proper closed convex cone $\gamma\subset \R^n$,
$\epsilon>0$ and an open neighborhood $W$ of  $x_0$ with $\xi_0\in
\Int(\gamma^{\circ})$ such that $(W+\gamma^a)\cap\ol{H_\epsilon}\subset X$
and
\begin{equation}\label{E:2}
H^j(X; \corps_{(x+{\gamma^a})\cap H_{\epsilon}}\otimes
F)=0\quad\mbox{for any $j\le k$ and $x\in W$.}
\end{equation}
\enum
\end{proposition}
\begin{proof}
We may assume $X=\R^n$.

\medskip
\noindent
(ii)${}_{k}\Rightarrow$ (i)${}_{k}$ is obvious.

\medskip
\noindent
(i)${}_k \Rightarrow$ (iv)${}_k$\quad
By the hypothesis, there exist distinguished triangles
$$ G\to F\to K\overset{+1}{\longrightarrow}\quad \mbox{and}\quad
G\to F'\to K'\overset{+1}{\longrightarrow}$$
in $D^b(\corps_X)$
such that $p\notin \Ss(K)$ and $p\notin \Ss(K')$.
By Lemma \ref{L:01},  there exist  an open neighborhood $W$ of
$x_0$,  a proper closed convex cone $\gamma$
such that $\xi_0\in \Int(\gamma^{\circ})$, and
$\epsilon>0$ such that
\eqn
&& H^j(X; \corps_{(x+\gamma^a)\cap H_{\epsilon}}\otimes K)
=H^j(X;\corps_{(x+\gamma^a)\cap H_{\epsilon}}\otimes K')
=0
\eneqn
for any $j\in\Z$ and $x\in W$.
Hence one has
\eqn
&&  
H^j(X; \corps_{(x+\gamma^a)\cap H_{\epsilon}}\otimes F)
\simeq H^j(X; \corps_{(x+{\gamma^a})\cap H_{\epsilon}}\otimes F').
\eneqn
Since 
$\corps_{(x+{\gamma^a})\cap H_{\epsilon}}\otimes F'$ belongs to 
$D^{>k} (\corps_X)$, 
we get \eqref{E:2}.

\medskip
\noindent
(i)${}_k \Rightarrow$ (iii)${}_{k,}{ }_{1}$\quad
Same proof as
(i)${}_k \Rightarrow$ (iv)${}_k$, replacing
$H^j(X; \corps_{(x+{\gamma^a})\cap H_{\epsilon}}\otimes G)$ 
with ${H}^j_{\{\phi\geq 0\}}(G)_x$
where $G= F'$, $F$, $K$, $K'$.

\medskip
\noindent
(iii)${}_{k,1} \Rightarrow $ (iii)${}_{k,\omega}$ is obvious.

\medskip
\noindent 
(iv)${}_k \Rightarrow$(ii)${}_k$\quad
To start with, note that (\ref{E:2}) entails \begin{equation}\label{E:5}
(\corps_{Z_{{\gamma}}}\circ F_{H_{\epsilon}})_W\in
D^{>k} (\corps_X).
\end{equation}

Let $\Delta$ denote the diagonal of $X\times X$.
Then the morphism
$\corps_{Z_{\gamma}}\to \corps_{\Delta}$
induces
the morphism in $D^{b}(\corps_X)$
\begin{equation}\label{E:4}
\corps_{Z_{{\gamma}}}\circ F_{H_{\epsilon}}\to  F_{H_{\epsilon}},
\end{equation}
which is an isomorphism in $D^{b}(\corps_X; p)$ 
by \cite[Theorem 7.1.2]{K-S1}.
Therefore, the composition
\eqn
&&(\corps_{Z_{\gamma}}\circ F_{H_{\epsilon}})_W
\to  (F_{H_{\epsilon}})_W\to F
\eneqn
is an isomorphism in $D^{b}(\corps_X;p)$
and  $(\corps_{Z_{\gamma}}\circ F_{H_{\epsilon}})_{W}$
belongs to $D^{>k}(\corps_X;p)$.

\medskip
\noindent (iii)${}_{k,\omega}$ $\Rightarrow$(iv)${}_k$\quad
We already know that (iv)${}_k$ is equivalent to
(i)${}_k$ for every $k$. Hence arguing by induction on $k$,
we may assume that (i)$_{k-1}$ holds.
Therefore we may assume $F\in D^{\ge k}(\corps_X)$.
Then we have
$$H^j(X; \corps_{(x+{\gamma^a})\cap H_{\epsilon}}\otimes
F)=0\quad\mbox{for any $j\leq k-1$}$$
and
\eqn
H^{k}(X; \corps_{(x+{\gamma^a})\cap H_{\epsilon}}\otimes F)
&\simeq& 
H^0(X; \corps_{(x+{\gamma^a})\cap H_{\epsilon}}\otimes H^{k}(F)),\\
H^{k}_{\{\varphi\ge0\}}(F)&\simeq&\Gamma_{\{\varphi\ge0\}}(H^{k}(F)).
\eneqn
We may assume that  
$\Int(\gamma)\neq \emptyset$,
${W}\times (\gamma^{\circ}\setminus\{0\}))\subset U$
and $(W+\gamma^a)\cap \overline{H_\epsilon}\subset W$.
Let 
$s\in \Gamma(X; \corps_{(x+{\gamma^a})\cap H_{\epsilon}}\otimes H^{k}(F))$.
Then there exists $y\in\R^n$ such that
$x+\gamma^a\subset y+\Int(\gamma^a)\subset W\cup(X\setminus H_\epsilon)$ and
$s$ extends to a section
\eqn
&&
\tilde s\in \Gamma\bigl(y+\Int(\gamma^a);\corps_{H_\epsilon}
\otimes{H}^{k}(F)\bigr)
\subset \Gamma\bigl(y+\Int(\gamma^a);{H}^{k}(F)\bigr).
\eneqn
Set $S=\supp(\tilde s)\subset H_\epsilon\cap(y+\Int(\gamma^a))$.
Then the following lemma asserts $S=\emptyset$,
and hence $H^{k}(X;\corps_{(x+{\gamma^a})\cap H_{\epsilon}}\otimes F)=0$.
\end{proof}

\begin{lemma}\label{lem:geo}
Let $\gamma$ be a proper closed convex cone in $\R^n$.
Let $\Omega$ be an open subset of $\R^n$ such that $\Omega+\gamma^a=\Omega$,
and let $S$ be a closed subset of $\Omega$ such that
$S\Subset\R^n$.
Assume the following condition:
for any $x\in \R^n$ and a real analytic function
$\varphi$ defined on $\R^n$,
the three conditions $S\cap\varphi^{-1}(\R_{<0})=\emptyset$,
$\varphi(x)=0$ and $d\varphi(x)\in\Int(\gamma^\circ)$
imply $x\notin S$.

Then $S$ is an empty set.
\end{lemma}
\begin{proof}
If $\gamma=\{0\}$, then by taking $\phi=0$, the lemma is trivially
true. Hence we may assume that $\{0\}\subsetneqq\gamma$.
Let us take $\xi$ such that 
$\gamma^\circ\setminus\{0\}\subset\{x;\lan x,\xi\ran>0\}$.
Then there is a real number $a$ such that $S\subset \{x;\lan x,\xi\ran>a\}$.
Set $H_-=\{x;\lan x,\xi\ran<a\}$.
By replacing $\Omega$ with $\Omega\cup H_-$,
we may assume from the beginning
$H_-\subset\Omega$.

For a proper closed convex cone $\gamma'$
such that $\gamma\setminus\{0\}\subset\Int(\gamma')$,
set $\Omega_{\gamma'}=\{x\in \Omega;x+\gamma'\in\Omega\}$ and
$S_{\gamma'}=S_{\gamma'}\cap \Omega_{\gamma'}$.
Since $S=\mathop\cup\limits_{\gamma'}S_{\gamma'}$, it is enough to show the 
assertion for $S_{\gamma'}$.
Since $\gamma'{}^\circ\setminus\{0\}\subset \Int(\gamma^\circ)$,
by replacing $\Omega$, $S$, $\gamma$ with
$\Omega'_{\gamma'}$, $S_{\gamma'}$ and $\gamma'$,
we may assume from the beginning
\eq
&&\Int(\gamma)\not=\emptyset,\\
&&\mbox{$S\cap\varphi^{-1}(\R_{<0})=\emptyset$,
$\varphi(x)=0$ and $d\varphi(x)\in\gamma^\circ\setminus\{0\}$}
\Longrightarrow x\notin S.
\eneq

Let us set 
$\psi(x)=\dist(x,\gamma^a):=\inf\{\Vert y-x\Vert\,;\,y\in\gamma^a\}$.
It is well known that
$\psi$ is a continuous function on $\R^n$,
and $C^1$ on $\R^n\setminus\gamma^a$.
More precisely for any $x\in \R^n\setminus\gamma^a$,
there exists a unique
$y\in\gamma^a$ such that
$\psi(x)=\Vert x-y\Vert$.
Moreover $d\psi(x)=\Vert x-y\Vert^{-1}(x-y)\in\gamma^\circ{}\setminus\{0\}$.
Furthermore $B_{\psi(x)}(y):=\{z\in\R^n;\Vert z-y\Vert<\psi(x)\}$
is contained in $\{z\in\R^n;\psi(z)<\psi(x)\}$.

For $\epsilon>0$, we set $\gamma^a_\epsilon=\{x\in\R^n;\psi(x)<\epsilon\}$.
Then $\gamma^a_\epsilon$ is an open convex set.
Moreover $\gamma^a_\eps+\gamma^a=\gamma^a_\eps$.
Set $\Omega_\epsilon=\{x;x+\gamma^a_\epsilon\subset\Omega\}$.
Then $\Omega=\cup_{\epsilon>0}\Omega_\epsilon$.
Set $S_\eps=S\cap\Omega_\eps$.
It is enough to show that $S_{\eps}=\emptyset$.

Assuming $S_{\eps}\not=\emptyset$, we shall derive a contradiction.
Let us take $x_0\in S_\eps$ and $v\in\Int(\gamma)$.
Set $V_t=x_0+\gamma_{\eps/2}^a+tv$ for $t\in\R$.
Then one has
\eq
&&\mbox{$V_t=\bigcup_{t'<t}V_{t'}$ and 
$\ol{V_t}=\bigcap_{t'>t}{V_{t'}}=\bigcap_{t'>t}\ol{V_{t'}}$,}\\
&&\mbox{$x_0\in V_t\cap S$ for $t\ge0$,
and $V_t\subset H_-$ for $t\ll0$,}\\
&&\mbox{$\ol{V_t}\subset\Omega$ for any $t\le0$,}\label{eq:clos}
\eneq
 Hence, for any compact set $K$ and $t\in\R$ such $K\cap\ol{V_t}=\emptyset$,
there exists $t'>t$ such that
$K\cap{V_{t'}}=\emptyset$.

Let us set $c=\sup\{t;V_t\cap S=\emptyset\}$.
Then $c\le 0$ and
$V_c\cap S=\emptyset$. 
By \eqref{eq:clos}, one has 
$\ol{V_c}\cap \ol{S}
\subset \ol{V_c}\cap S$.        
Since $\ol{S}$ is a compact set, there exists
$x_1\in S\cap\partial V_c$.
Here $\partial V_c:=\ol{V_c}\setminus V_c$ is the boundary of
$V_c$.
As seen before, there exists
a ball
$B_{\eps/2}(y):=\{x;\Vert x-y\Vert<\eps/2\}$ such that
$B_{\eps/2}(y)\subset V_c$, $\Vert x_1-y\Vert=\eps/2$
and $x_1-y\in\gamma^\circ$.
This is a contradiction by taking $\phi(x)=\Vert x-y\Vert^2-(\eps/2)^2$.
\end{proof}

\begin{definition}\label{D:1}
\bnum
\item 
Let $F\in D^b(\corps_X)$.
The closed conic subset $\Ss_k(F)$ of $T^*X$ is defined by:
$p\notin \Ss_k(F)$ if and only if 
$F$ satisfies the equivalent conditions in
Proposition \ref{L:1}.
\item
Let $p\in T^*X$ and $k\in\Z$. 
Then $ D^{>k}(\corps_X; p)$ denotes the full
additive subcategory of $ D^{b}(\corps_X; p)$ 
consisting of $F$ satisfying $p\notin\Ss_k(F)$.
We write sometimes $D^{\ge k}(\corps_X; p)$ for $D^{>k-1}(\corps_X; p)$.
\enum
\end{definition}
Note that  $\Ss_k(F)\cap T^*_XX=\pi(\Ss_k(F))=\Supp(\tau^{\le k}F)$.

\begin{remark}
The truncated microsupport has the following properties,
similarly to those of the  microsupport.
\bnum
\item For any $F\in D^b(\corps_X)$, one has $\Ss_k(F[n])=\Ss_{k+n}(F)$.
\item If $F'\to F\to F''\xrightarrow{+1}$ is a distinguished triangle, then one has
\eq\label{eq:triangle}
&&\ba{l}
\Ss_k(F)\subset \Ss_k(F')\cup \Ss_k(F''),\\[3pt]
\bigl(\Ss_k(F')\setminus \Ss_{k-1}(F'')\bigr)
\cup \bigl(\Ss_k(F'')\setminus \Ss_{k+1}(F')\bigr)
\subset \Ss_k(F).
\ea
\eneq
\enum
\end{remark}

\begin{remark}
\bnum
\item
If $F\in D^{>k}(\corps_X)$, then
$\Ss_k(F)=\emptyset$.
\item
If $F\in D^{\le k}(\corps_X)$, then
$\Ss_{k+d_X}(F)=\Ss(F)$.
Here $d_X$ is the dimension of $X$.
\enum
The last statement follows from 
the characterization (iv)$_k$ in Proposition \ref{L:1} and
the fact that $H^j(X;F)$ vanishes for
any $F\in D^{\le k}(\corps_X)$ and $j>k+d_X$.
\end{remark}

\begin{remark}
It is not true that $F\in D^{\ge k}(\corps_X;p)$
implies the existence of a morphism $F\to F'$ in $D^b(\corps_X)$
which is an isomorphism in $D^b(\corps_X;p)$
and $F'\in D^{\ge k}(\corps_X)$.
For example take $X=\R$, $p=(0\,;1)$, $Z=\{x\in X;x<0\}$,
$F_1=\corps_{\{0\}}$ and $F=\corps_Z[1]$.
Then there is a morphism $F_1\to F$ which is an isomorphism
in $D(\corps_X;p)$. Hence one has $F\in D^{\ge0}(\corps_X;p)$.
Assume that there is 
a morphism $u\colon F\to F'$ in $D^b(\corps_X)$
which is an isomorphism in $D^b(\corps_X;p)$
and $F'\in D^{\ge0}(\corps_X)$.
Since $H^0(F)=0$, the morphism  $H^0(F_1)\to H^0(F')$
vanishes, and hence the composition
$F_1\to F{\xrightarrow{\,u\,}}F'$ vanishes.
This is a contradiction.
\end{remark}

\begin{examples}
\bnum
\item 
One has 
\eqn
\Ss_k(\corps_X)=\left\{ \begin{array}{ll}
\emptyset &\mbox{for $k<0$,}\\[3pt]
T^*_XX &\mbox{for $k\geq 0$.}
\end{array}\right.
\eneqn
\item
Let $X=\R$ and $Z_1=\{x\in X; x\ge0\}$, $Z_2=\{x\in X; x>0\}$.
Then one has
\eqn
\ba{l}
\Ss_k(\corps_{Z_1})=
\left\{\ba{ll}
\emptyset &\mbox{for $k<0$,}\\[5pt]
\{(x;\xi);\xi=0,\,x\ge0\}\\
\hspace*{50pt}\cup\{(x;\xi);x=0,\,\xi\ge0\}\quad
&\mbox{for $k\geq0$,}
\end{array}\right.\\[30pt]
\Ss_k(\corps_{Z_2})=
\left\{\ba{ll}
\emptyset &\mbox{for $k<0$,}\\[5pt]
\{(x;\xi);\xi=0,\,x\ge0\}&\mbox{for $k=0$,}\\[5pt]
\{(x;\xi);\xi=0,\,x\ge0\}\\
\hspace*{50pt}\cup\{(x;\xi);x=0,\,\xi\le0\}\quad
&\mbox{for $k\geq 1$.}
\end{array}\right.
\ea
\eneqn
\item
Let $X$ be a complex manifold. Then 
\eqn
\Ss_k(\sho_X)=\left\{ \begin{array}{ll}
\emptyset &\mbox{for $k<0$,}\\[3pt]
T_X^*X&\mbox{for $k=0$,}\\[3pt]
T^*X &\mbox{for $k\geq 1$.}
\end{array}\right.
\eneqn
\item
Let $M$ be a real analytic manifold, $X$ a complexification of $M$,
$\shm$ a coherent $\shd_X$-module, and let $\shb_M$ denote the sheaf of
Sato's hyperfunctions on $M$. 
Regarding $T^*M$ as a subset of $T^*X$, one has
$$SS_0\bigl(\rhom[\shd_X](\shm,\shb_M)\bigr)
\subset \car(\shm)\cap T^*M.$$
This follows immediately from the Holmgren theorem.
\enum
\end{examples}

\section{Functorial properties}\label{section:3}

In this section, we  study in Propositions \ref{P:3}--\ref{P:5} below,   
the behavior of $\Ss_k$ under external tensor products, 
proper direct image and
smooth inverse image. These properties are proved similarly to   
the corresponding properties of the  microsupport 
(cf.\ Chapter V of \cite{K-S1}), 
using Proposition \ref{L:1}.

\begin{proposition}\label{P:3}
Let $X$ and $Y$ be real analytic manifolds.
Then for $F\in D^{b}(\corps_X)$, $G\in D^{b}(\corps_Y)$ and
$k\in\Z$, one has
\begin{equation}\label{E:09}
\Ss_k(F\etens G)\subset 
\bigcup_{i+j=k}\Ss_i(F)\times \Ss_j(G).  
\end{equation}
\end{proposition}
\begin{proof}
Let us show that $(p,p')\notin
\bigcup_{i+j=k}\Ss_i(F)\times \Ss_j(G)$
implies
$(p,p')\notin \Ss_k(F\etens G)$.
Since
$\Ss_k(F\etens G)
\subset \Ss(F\etens G)
\subset \Ss(F)\times \Ss(G)$,
we may assume that
$(p,p')\in  \Ss(F)\times \Ss(G)$.
Since $\Ss_i(F)=\Ss(F)$ for $i\gg0$ and $\Ss_i(F)=\emptyset$ for $i\ll0$, 
there exists $i$ such that $p\in \Ss_i(F)$ and $p\notin \Ss_{i-1}(F)$.
Set $j=k-i$.
Then $p'\notin \Ss_j(G)$.
Hence there exist a morphism $F'\to F$ 
which is an isomorphism in $D^b(\corps_X;p)$
with $F'\in D^{>i-1}(\corps_X)$, and
and $G'\to G$ which is an isomorphism in $D^b(\corps_Y;p')$
with $G'\in D^{>j}(\corps_X)$.
Hence $F'\etens G'
\to F\etens G$
is an isomorphism in $D^b(\corps_{X\times Y};(p,p'))$ and
$F'\etens G'\in D^{>i+j}(\corps_{X\times Y})$.
\end{proof}

\begin{proposition}\label{P:4}
Let $f\colon X\to Y$ be a morphism of real analytic manifolds
and let $F\in D^{b}(\corps_X)$ such that $f$ is proper on the support of $F$.
Then for any $k\in\Z$,
\begin{equation}\label{E:10}
\Ss_k (\roim{f}(F))\subset f_{\pi}\opb{{f_d}}(\Ss_k(F)).
\end{equation}
The equality holds in case f is a closed embedding.
\end{proposition}

\begin{proof}
We shall follow the method of proof of Proposition 5.4.4. of \cite {K-S1}.

Let $y\in Y$ and let $\phi$ be a real $C^1$-function on $Y$
such that $\phi(y)=0$ and
$d(\phi\circ f)(x)\notin \Ss_k(F)$ for every $x\in \opb{f}(y)$.
Therefore
\eqn
&&
H^j\rsect_{\{\phi\circ f\geq 0\}}(F)\vert_{\opb{f}(y)}=0\quad
\mbox{for any $j\le k$.}
\eneqn
We have
\eqn
H^j \rsect_{\{\phi\geq 0\}}(\roim{f}(F))_y
&\simeq&H^j\roim{f}(\rsect_{\{\phi\circ f\geq0\}}(F))_y\\
&\simeq& H^j(\opb{f}(y);\rsect_{\{\phi\circ f\geq 0\}}(F))=0
\eneqn
for every $j\le k$. This proves (\ref{E:10}).

Let us now assume that $f$ is a closed embedding. Let $p\notin \Ss_k(Rf_*F)$.
 We may assume that  $Y$ is a real vector space
and $X$ is a linear subspace of $Y$.
Let $\gamma\subset Y$, $W\subset Y$,  $\epsilon$  be chosen as in
Proposition \ref{L:1} (iv)${}_k$
with respect to $p$ and $Rf_*F$, that is, $$ H^j(Y;
\corps_{(x+{\gamma^a})\cap H_{\epsilon}}\tens \roim{f}F)=0
\quad\mbox{for any $j\le k$ and $x\in W$.}$$
Since $f$ is a closed embedding, one has
\eqn
H^j(Y;\corps_{(x+{\gamma^a})\cap H_{\epsilon}}\tens Rf_*F)
&\simeq&
 H^j(X; \corps_{(x+\gamma^a)\cap H_{\epsilon}\cap X}\tens F).
\eneqn
Hence one has
\eqn
&&
H^j(X;\corps_{(x+{\gamma^a}\cap X)\cap (H_{\epsilon}\cap X)}\tens F)=0
\quad\mbox{for any $j\le k$,}
\eneqn
and the interior of 
the polar set of $\gamma\cap X$ contains $f_d \opb{f_\pi}(p)$,
and therefore
$\Ss_k(F)\cap f_d \opb{f_\pi}(p)=\emptyset$.
\end{proof}

\begin{proposition}\label{P:5}
Let $X$ and $Y$ be real analytic manifolds and let $f\colon X\to Y$
be a smooth morphism. Let $G\in D^b(\corps_Y)$. Then, for any $k\in\Z$,
\begin{equation}\label{E:11}
\Ss_k (\opb{f}G)= f_d\opb{f_\pi}(\Ss_k (G)).
\end{equation}
\end{proposition}
\begin{proof}
The problem being
local on $X$, we may assume that $X=Y\times Z$,
$Y$ and $Z$ are vector spaces and $f$ is the projection. Then we have to
show 
\eq\label{eq:ssopb}
&& \Ss_k(G\etens \corps_Z)=\Ss_k(G)\times T^*_ZZ.
\eneq
The inclusion  $\subset$ in \eqref{eq:ssopb}
is a particular case of Proposition \ref{P:3}. Let us prove the converse
inclusion. 
Since $\Ss_k(G\etens \corps_Z)\subset \Ss(G)\times T^*_ZZ$, 
it is enough to show that
$(p,p')\notin\Ss_k(G\etens \corps_Z)$ ($p'=(z_0;0)$)
implies $p\notin\Ss_k(G)$.

Setting $p=(x_0;\xi_0)$,
there exist a proper closed convex cone $\gamma$ in $X\times Z$
and an open neighborhood $W$ of $x_0$ and $\epsilon$ 
such that $(\xi_0,0)\in\Int(\gamma^\circ)$ and
$$H^j(X\times Z;\corps_{(H_\eps\times Z)\cap((x,z_0)+\gamma^a)}
\tens(G\etens \corps_Z))=0
\quad\mbox{for $j\le k$ and $x\in W$.}$$
One has
\eqn
&&\rsect\Bigl(X\times Z;\corps_{(H_\eps\times Z)\cap((x,z_0)+\gamma^a)}
\tens(G\etens \corps_Z)\Bigr)
\\
&&
\hspace*{80pt}
\simeq
\rsect\Bigl(X;\roim{f}\bigl(\corps_{(H_\eps\times Z)\cap((x,z_0)+\gamma^a)}
\tens(G\etens \corps_Z)\bigr)\Bigr)
\eneqn
and 
$\roim{f}\bigl(\corps_{(H_\eps\times Z)\cap((x,z_0)+\gamma^a)}
\tens(G\etens \corps_Z)\bigr)
\simeq\roim{f}(\corps_{(x,z_0)+\gamma^a})\tens\corps_{H_\eps}\tens G$.
Since one has
$$\roim{f}(\corps_{(x,z_0)+\gamma^a})\simeq \corps_{x+f(\gamma^a)},$$
we obtain
$$H^j(X;\corps_{H_\eps\cap(x+f(\gamma)^a)}
\tens G)=0
\quad\mbox{for $j\le k$ and $x\in W$.}$$
Then the assertion follows from the fact that
$f(\gamma)$ is a proper closed convex cone such that
$\xi_0\in\Int(f(\gamma)^\circ)$.
\end{proof}

\section{Estimates for the truncated microsupport}
\label{S:4}

Let $Y$ be  a  smooth submanifold of $X$.
In this section  we will give an estimate for $\Ss_k(F)\cap \opb{\pi}(Y)$.
Recall that $\mu_Y (F)$ denotes the microlocalization of $F$ along $Y$.
Note that for $F\in D^{\ge k}(\corps_X)$,
$H^k(\mu_Y (F))\simeq H^0\bigl(\mu_Y (H^k(F))\bigr)$ 
is a subsheaf of $\pi^{-1}H^k(F)\vert_{T^*_YX}$ and
\eq\label{eq:Sato}
&&H^k(\mu_Y (F))_p\simeq
\ \bigl\{s\in H^k(F)_{\pi(p)}\,;\,
C_Y(\supp(s))_{\pi(p)}\setminus\{0\}\\[1pt]
&&\hspace*{180pt}\subset\{v\in (T_YX)_{\pi(p)};\lan v,p\ran>0\}\bigr\}\nonumber
\eneq
for $p\in T^*_YX$.

The following result is a generalization of
\cite[Theorem 5.7.1]{K-S2} to $\Ss_k$.
\begin{theorem}\label{T:5}
Let $X$ be a real analytic manifold and $Y$ a smooth submanifold. 
Let $F\in D^b (\corps_X)$. Then
\eq\label{E:16}
&&\Ss_k (F)\cap T^*_Y X=
\bigl(T^*_YX\cap\overline{\Ss_k (F)\setminus \opb{\pi}(Y)}\bigr)
\cup \Supp (\tau^{\le k}\mu_Y (F)).
\eneq
\end{theorem}

\begin{proof}
It is evident that the right hand side of \eqref{E:16}
is contained in the left hand side.
Let us show the converse inclusion.
Assuming that  $p\in T^*_Y X$ satisfies
\eqn
&& p\notin \overline{\Ss_k (F)\setminus \opb{\pi}(Y)}
\cup \Supp\bigl(\tau^{\le k}(\mu_Y (F))\bigr),
\eneqn
we shall prove $p\notin \Ss_k(F)$.
Arguing by induction on $k$, one has
$p\notin \Ss_{k-1}(F)$ and
by Proposition \ref {L:1} (ii)${}_{k-1}$  we may assume that
$F\in D^{\geq k}(\corps_X)$.

There exists an open conic
neighborhood $U$ of $p$ in $T^*X$
such that $U\cap \Ss_k (F) \subset \pi^{-1}(Y)$ and
${H}^j (\mu_Y (F))\vert_U =0$ for any $j\le k$.
Furthermore, we may assume that $X=\R^n$,
$Y$ is a linear subspace of $X$,
$p=(x_0;\xi_0)$.
Let us take 
an open neighborhood $W$ of
$x_0$, a proper closed convex cone $\gamma$ and $\eps>0$ such that
$W\times \Int(\gamma^{\circ})\subset U$,
$\xi_0\in\Int(\gamma^{\circ})$
and $(W+\gamma^a)\cap\ol{H_\eps}\subset W$.
Hence one has
\eqn
&& H^j(X;\corps_{{(x+{\gamma^a})}\cap H_{\epsilon}}\otimes F)
\simeq H^{j-k}(X;\corps_{{(x+{\gamma^a})}\cap H_{\epsilon}}\otimes H^k(F))
\quad\mbox{for any $j\leq k$.}
\eneqn
Thus it is enough to check that
\eq\label{E:016}
\Gamma\bigl(X;\corps_{(x+{\gamma^a})\cap H_{\epsilon}}\otimes H^{k}(F)
\bigr)=0.
\eneq

Let 
$s\in \Gamma(X; \corps_{(x+{\gamma^a})\cap H_{\epsilon}}\otimes H^{k}(F))$.
Then there exists an open set $\Omega_0$ such that
$\Omega_0+\gamma^a=\Omega_0$, $x+{\gamma^a}\subset\Omega_0$
and $s$ extends to a section
$\tilde s\in\Gamma(\Omega_0;\corps_{H_{\epsilon}}\otimes H^{k}(F))$.
Moreover we may assume that
$(\Omega_0\cap \ol{H_\eps})\times\gamma^\circ\subset U$.
Set $S=\supp(\tilde s)$.
Then $S\setminus (Y+\gamma)$ satisfies the condition
in Lemma \ref{lem:geo} with $\Omega=\Omega_0\setminus (Y+\gamma)$.
Hence we have
$S\setminus (Y+\gamma)=\emptyset$ and hence $S\subset Y+\gamma$.
Since $H^{k}(\mu_Y(F))\vert_U=0$ and
$$C_Y(Y+\gamma)\subset Y\times(Y+\gamma)
\subset\{v\in T_YX;\lan v,\xi_0\ran>0\}\cup (Y\times\{0\}),$$
 the formula \eqref{eq:Sato}
implies $\tilde s\vert_Y=0$.
One has therefore $S\cap Y=\emptyset$.
Then $S$ satisfies the condition in Lemma \ref{lem:geo},
and we can conclude $S=\emptyset$,
which implies $s=0$.
\end{proof}

We shall
need the following definition:

\begin{definition}\label{D:4}
Let $Y$ be a closed submanifold of $X$, let $k\in\Z$ and 
let $F\in D^b (\corps_X)$. The closed subset $\Ss_k^Y(F)$ of $\opb{\pi}(Y)$
is defined by: $p\notin \Ss_k^Y(F)$ if and only if 
there exists an open conic neighborhood
$U$ of $p$ in $\opb{\pi}(Y)$ satisfying the following two conditions:
\bnum
\item 
$\tau^{\le k}\mu _{Y}(F)\vert_{U\cap T^*_YX}=0$,
\item 
for any smooth real analytic hypersurface
$Z$ of $Y$,
$$\tau^{\le k}\mu _{Z}(F)\vert_{U\cap T^*_ZX\setminus T^*_YX}=0.$$
\enum
\end{definition}
We remark that
$\Ss_k^Y(F)$ is a conic closed set obviously contained in
$\Ss_k(F)\cap \opb{\pi}(Y)$.

\begin{theorem}\label{T:6}
Let $X$ be a real analytic manifold and $Y$ a closed submanifold. 
Let $F\in D^b (\corps_X)$. Then
\eq\label{E:14}
\Ss_k (F)=\overline{\Ss_k (F)\setminus \opb{\pi}(Y)}\cup \Ss_k^Y(F).
\eneq
\end{theorem}
\begin{proof}
The left hand side obviously contains the right hand side.
Let us prove the converse inclusion.
Assuming that $p\in \opb{\pi}(Y)$ satisfies
$p\notin \overline{\Ss_k (F)\setminus \opb{\pi}(Y)}\cup \Ss_k^{Y}(F)$,
let us show $p\notin \Ss_k(F)$.

If $p\in T^*_YX$, Theorem \ref{T:5} implies the assertion.
Hence we may assume $p\notin T^*_YX$.
Let $U$ be an open conic neighborhood of $p$ in $T^*X$ such that
$\Ss_k(F)\cap U\subset \pi^{-1}(Y)$ and $U\cap \Ss_k^Y(F)=\emptyset$. 

We may assume that $X=\{x=(u,v,t);u\in \R^n,\,v\in \R^m,\,t\in \R\}$,
$Y=\{(u,v,t)\in X;u=0\}$ and $p=((0,0,0);(0,0,1))$.
We may assume $W\times\gamma^\circ\subset U$
with $\gamma=\{t\ge \sqrt{\Vert u\Vert^2+\Vert v\Vert^2}\}$
and an open neighborhood $W$ of the origin.
Set $H_\epsilon=\{(u,v,t);t>-\epsilon\}$,
and choose $W$ and a sufficiently small $\epsilon$ such that
$(x+\gamma^a)\cap \ol{H_\epsilon}\subset W$ for any $x\in W$.

By the induction on $k$, we may assume $F\in D^{\ge k}(\corps_X)$.
It is enough to show for any $x_0\in W$
$$\Gamma\bigl(x_0+\Int(\gamma^a);\corps_{H_\epsilon}\otimes H^{k}(F)\bigr)=0.$$
Let $s\in \Gamma\bigl(x_0+\Int(\gamma^a);
\corps_{H_\epsilon}\otimes H^{k}(F)\bigr)$.
Set $S=\supp(s)\subset\bigl(x_0+\Int(\gamma^a)\bigr)\cap H_\epsilon$.
Assuming $S\not=\emptyset$, we shall derive a contradiction.
Set $x_0=(u_0,v_0,t_0)$
and set $\phi(x)=\Vert u-u_0\Vert^2+\Vert v-v_0\Vert^2-(t-t_0)^2$.
Then $\Omega:=x_0+\Int(\gamma^a)=\{x;t<t_0, \phi(x)<0\}$,
and $\varphi(\Omega\cup H_\eps)$ is bounded from below.
Moreover one has $d\phi(x)\in\Int(\gamma^\circ)$ for any $x\in\Omega$,
and $d\phi(x)\notin T^*_YX$ for any $x\in\Omega\cap Y$.
Set $c=\inf\{\phi(x);x\in S\}<0$.
Since $\phi\vert_S\colon S\to\R_{<0}$ is a proper map,
one has $c\in \phi(S)$.
Let $Z_c$ be the closed subset
$\{x=(0,v,t)\in Y;t<t_0,\,\phi(x)=c\}$ of $Y$ and
set
$\Omega'=\Omega\setminus ({Z_c}+\gamma)$.
Since one has
$${Z_c}+\gamma
=\bigl\{(t,u,v);
t-t_0\ge \Vert u\Vert-\sqrt{\Vert v-v_0\Vert^2+\Vert u_0\Vert^2+c}\bigr\},$$
$\phi$ takes values smaller than $c$ on $Y\setminus ({Z_c}+\gamma)$,
and hence
$S':=S\cap\Omega'$ does not intersect $Y$.
Therefore $S'$ satisfies the condition
in Lemma \ref{lem:geo}.
Hence $S'=\emptyset$, which means 
$S\subset {Z_c}+\gamma$.
Since $C_{Z_c}({Z_c}+\gamma)_x\setminus\{0\}\subset\{d\phi(x)>0\}$
for any $x\in Z_c$,
$H^0(\mu_{Z_c}H^{k}(F))\vert_{U\setminus T^*_YX}=0$ implies
$s\vert_{Z_c}=0$. Therefore $S\cap{Z_c}=\emptyset$.
Since $\{x\in({Z_c}+\gamma)\cap\Omega;\phi(x)=c\}\subset{Z_c}$,
one has $c\notin\phi(S)$, which is a contradiction. 
\end{proof}

\section{Applications to $\shd$ -Modules}\label{S:5}
In this section, $X$ denotes a complex manifold.

Before stating our main result, let us recall a classical lemma on the
vanishing of the microlocalization of $\sho_X$ along submanifolds.

\begin{lemma}\label{lem:van}
Let $Y$ be a closed complex submanifold of codimension $d$ of $X$  and let 
$S$ be a smooth real analytic hypersurface of $Y$.
Then
\eq
&&H^k(\mu_Y(\sho_X))=0\quad\mbox{for any $k\neq d$},\label{eq:skk}\\
&&H^k(\mu_S(\sho_X))\vert_{T^*_SX\setminus T^*_YX}=0
\quad\mbox{for any $k\le d$.}
\label{eq:ks2}
\eneq
\end{lemma}
\begin{proof}
The vanishing property
\eqref{eq:skk} is proved in \cite{S-K-K} and (by a different
method) in \cite{K-S2}, Proposition 11.3.4.

\medskip
The vanishing property \eqref{eq:ks2} follows from 
\cite{K-S2}, Proposition 11.3.1.
Let us recall this statement. Let $p\in T^*_SX$. Set 
$E_p=T_p(T^*X)$, $\lambda_S=T_p(T^*_SX)$,
$\lambda_0=T_p(\opb{\pi}\pi(p))$,
and denote by $\nu$ the
complex line in $E_p$, the tangent space to the Euler vector field in $T^*X$
at $p$. Let $c$ be the real codimension of the real submanifold $S$
 and let $\delta$ denote the complex dimension of
$\lambda_S\cap \sqrt{-1}\lambda_S\cap\lambda_0$. 

The result of
loc.\ cit.\ asserts that if the real dimension of
$\lambda_S\cap\nu$ is $1$, then
\eqn
&& H^j(\mu_S(\sho_X))_p=0\quad\mbox{for $j<c-\delta$.}
\eneqn 
(The result
in loc.\ cit.\ is more precise, involving the signature of
the Levi form.)
If $p\in T^*_SX\setminus T^*_YX$,
the real dimension of $\lambda_S\cap\nu$ is $1$.
Since $c=2d+1$ and $\delta=d$, we get the desired result.
\end{proof}

Now we are ready to prove the following proposition.

\begin{proposition}\label{T:7}
Let $X$ be a complex manifold, let $\shm$ be
a coherent $\shd_X$-module and
let $S$ be a closed complex analytic subset of $X$ 
with $\codim_XS\ge d$. 
Set
$F=R\hom[\shd_X](\shm,\sho_X)$. Then
\bnum
\item
$\Ss_{d-1}(F)=\ol{\Ss_{d-1}(F)\setminus\opb{\pi}(S)}$.
\item
Let $S'$ be a closed complex analytic subset of $S$
such that $\codim_XS'>d$ and $S_0:=S\setminus S'$
is a non-singular subvariety of codimension $d$.
Then
\eqn
&&\Ss_{d}(F)=\ol{\Ss_{d}(F)\setminus\opb{\pi}(S)}
\cup\ol{\Supp\bigl(\tau^{\leq d}\mu_{S_0}(F\vert_{X\setminus S'})\bigr)}.
\eneqn
In particular one has
$$\Ss_{d}(F)=\ol{\Ss_{d}(F)\setminus\opb{\pi}(S)}\cup
\ol{\Ss_d(F)\cap T^*_{S_0}X}.$$
\enum
\end{proposition}
\begin{proof}
(i) By the induction on the codimension of $S$,
we may assume that $S$ is non-singular.
By Theorem \ref{T:6} one has 
\eqn
&& \Ss_{d-1}(F)=\ol{\Ss_{d-1}(F)\setminus\opb{\pi}(S)}\cup \Ss_{d-1}^S(F).
\eneqn
Hence it is enough to show that $\Ss_{d}^S(F)=\emptyset$, or
equivalently
$H^j(\mu_S(\sho_X))=0$ for $j<d$
and 
$H^j(\mu_Z(\sho_X))\vert_{T^*_ZX\setminus T^*_SX}=0$ for $j<d$
for any real analytic hypersurface $Z$ of $S$.
This is a consequence of Lemma \ref{lem:van}.

\medskip
\noindent
(ii) By (i), we may assume that $S$ is non-singular of codimension $d$.
Hence it is enough to show
$\Ss_{d}^S(F)=\Supp\bigl(\tau^{\leq d}\mu_{S}(F)\bigr)$.
By the definition, we are reduced to proving
$H^j(\mu_Z(\sho_X))\vert_{T^*_ZX\setminus T^*_SX}=0$ for $j\le d$ 
and for any real analytic hypersurface $Z$ of $S$.
This is again a consequence of Lemma \ref{lem:van}.
\end{proof}

\begin{remark}
When $S$ is a closed smooth hypersurface, the inclusion 
\eqn
&&\Ss_{1}(F)\subset\ol{\Ss(F)\setminus\opb{\pi}(S)}
\cup(\Ss(F)\cap T^*_SX)
\eneqn
was obtained in \cite{To}.
\end{remark}

Let $\Omega$ be an open subset of $X$. We shall say for short that
$\Omega$ has a smooth boundary $\partial \Omega$ 
if there exists a real C$^1$-function
$\phi$ such that $d\phi\neq 0$ on the set $\{\phi=0\}$ and
$\Omega=\{x\in X;\phi(x)<0\}$. 

\begin{corollary}\label{co:7}
Let $\Omega$ be an open subset of $X$ with smooth boundary, let 
$\shm$, $F$, $S$ and $S_0$ be as in Proposition \ref{T:7} and let
$\Lambda$ be a closed conic subset of $T^*X$.
Assume that
\eqn
&&\car(\shm)\subset\Lambda\cup\opb{\pi}(S),\\
&& \ol{T^*_{S_0}X}\cap T^*_{\partial \Omega}X\subset T^*_XX,\\
&&\Lambda\cap T^*_{\partial \Omega}X\subset T^*_XX.
\eneqn
Then one has
\eqn
&& \Ss_{d}(F)\cap T^*_{\partial \Omega}X\subset T^*_XX.
\eneqn
In particular one has
\eqn
&& H^j(\rsect_{X\setminus\Omega}(F))\vert_{\partial \Omega}=0
\quad\mbox{for $j\leq d$.}
\eneqn
\end{corollary}

\begin{example}
Under the situation of Corollary \ref{co:7},
assume further that $\Omega$ is pseudo-convex.
Let us denote by
$j\colon\Omega\hookrightarrow X$ the open embedding. Then 
$H^k(\roim{j}j^{-1}\sho_X)=0$ for $k\neq 0$, and 
$\rsect_{X\setminus\Omega}(\sho_X))\vert_{\partial \Omega}[1]$ is
concentrated in degree $0$.
Let us set for short:
\eqn
\sho_X^+/\sho_X&=&(\oim{j}j^{-1}\sho_X/\sho_X)\vert_{\partial \Omega}\\
&\simeq&\rsect_{X\setminus\Omega}(\sho_X))\vert_{\partial \Omega}[1].
\eneqn
Applying Corollary \ref{co:7}, we find that
\eqn
&&\ext[{\shd_X}]^j(\shm,\sho_X^+/\sho_X)=0\quad\mbox{for $j<d$.}
\eneqn
\end{example}

\begin{example}
Let $P$ be a differential operator on $X$
whose principal symbol $\sigma(P)$
has the form $a(x)q(x,\xi)$ with $a\in \sho_X(X)$ and
$q\in\sho_{T^*X}(T^*X)$.
Then taking  $\shd_X/\shd_X P$ as $\shm$,
the solution complex $F$ is
$\sho_X{\xrightarrow{P}}\sho_X$, where $\sho_X$ is at degree $0$ and $1$.
Taking $a^{-1}(0)$ and $q^{-1}(0)$ as $S$ and $\Lambda$,
Corollary \ref{co:7} implies
$$\Ss_1(F)\subset q^{-1}(0)\cup \ol{\{(x;\xi); a(x)=0, \xi\in\C da(x)\}}.$$
By \eqref{eq:triangle} and the distinguished triangle
$$\Ker(\sho_X\xrightarrow{P}\sho_X)\to 
F\to\Coker(\sho_X\xrightarrow{P}\sho_X)[-1]\xrightarrow{+1},$$
one has also
$$\Ss_1\bl(\Ker(\sho_X\xrightarrow{P}\sho_X)\br)
\subset q^{-1}(0)\cup \ol{\{(x;\xi); a(x)=0, \xi\in\C da(x)\}}.$$
\end{example}

\medskip
Finally one has the following theorem
which calculates $\Ss_k(F)$ in terms of $\car(\shm)$.
Here, for a closed complex subset $Z$ of $X$,
$T^*_ZX$ means $\ol{T^*_{Z_\reg}X}$ where
$Z_\reg$ is the non-singular locus of $Z$.

\begin{theorem}\label{th:sskFcharM}
Let $\shm$ be a coherent $\shd_X$-module,
and let $F$ be the solution complex $\rhom[{\shd_X}](\shm,\sho_X)$.
Let $\car(\shm)=\bigcup_{\alpha\in A}V_\alpha$
be the decomposition of $\car(\shm)$ into irreducible components.
Let $Y_\alpha$ be the irreducible complex analytic subset
$\pi(V_\alpha)$ of $X$.
Then for any integer $k$ one has
\eq\label{eq:ss}
&&\Ss_k(F)=\Bigl(\mathop\cup\limits_{\codim Y_\alpha<k}V_\alpha\Bigr)
\cup\Bigl(\mathop\cup\limits_{\codim Y_\alpha=k}T^*_{Y_\alpha}X\Bigr).
\eneq
\end{theorem}

\begin{proof}
The inclusion $\,\subset\,$ is a consequence of Proposition \ref{T:7}.
Let us show the converse inclusion.
Note that both sides are empty sets for $k<0$.
Hence arguing by induction on $k$,
we can assume that \eqref{eq:ss} holds for $k-1$.
Hence it is enough to show:
\begin{subequations}\label{eq:main}
\begin{gather}
\mbox{if $\codim Y_\alpha=k-1$ and $V_\alpha\not=T^*_{Y_\alpha}X$, then
$V_\alpha\subset \Ss_k(F)$,}\tag{\ref{eq:main}.i}\label{eq:a}\\
\mbox{if $\codim Y_\alpha=k$, then $T^*_{Y_\alpha}X\subset \Ss_k(F)$.}
\tag{\ref{eq:main}.ii}\label{eq:b}
\end{gather}
\end{subequations}

In both cases, we may assume that $Y:=Y_\alpha$ 
is a non-singular subvariety.
Let $j\cl Y\to X$ be the inclusion map.

\medskip
\noindent
{\em Proof of \eqref{eq:a}}\quad
It is enough to show
that for any open subset $U$ of $T^*X$ with a non empty intersection with
$V_\alpha$, $\Ss_k(F)\cap U$ is non empty.
We may assume that
$\car(\shm)\cap U=V_\alpha\cap U$
and $V_\alpha\cap U\to Y$ is a smooth morphism.
Since $V_\alpha\subset\pi^{-1}(Y)$,
we may assume, by shrinking $U$ if necessary, that
$\she_X^\R\tens[\shd_X]\shm\vert_U\simeq \she_{X\hot Y}^\R
\tens[\she_Y]\shn\vert_U$
for a coherent $\she_Y$-module $\shn$ (\cite{S-K-K}).
For any smooth complex hypersurface $Z$ of $Y$, one has by \cite{S-K-K}
\eqn
H^{k}(\mu_Z(F))&\simeq&
\hom[\shd_X]\bl(\shm,H^k(\mu_Z(\sho_X))\br)\\
&\simeq&\hom[\she_X^\R]\bl(\she_X^\R\tens[\shd_X]\shm,H^k(\mu_Z(\sho_X))\br)\\
&\simeq&
\hom[\she_X^\R]\bl(\she_{X\hot Y}^\R\tens[\she_Y]\shn,H^k(\mu_Z(\sho_X))\br)\\
&\simeq&j_d{}_*j_\pi^{-1}\hom[\she_Y]\bl(\shn, H^1(\mu_Z(\sho_Y))\br)
\eneqn
on $U\cap T^*_ZX$.
Hence the result follows from the following
lemma which is an easy consequence of
the classification theorem for coherent $\she$-modules
at generic points of their supports.
\begin{lemma}
Let $V$ be a non-empty {\rm(}locally closed\,{\rm)}
smooth submanifold of $\dot{T}^*Y$
such that $V\to Y$ is smooth,
and let $\shn$ be a coherent $\she_Y$-module 
defined on a neighborhood of $V$
such that $\Supp(\shn)=V$.
Then there is a  smooth complex hypersurface $Z$ of $Y$
such that
$\hom[\she_Y](\shn,H^1(\mu_Z(\sho_Y)))\vert_{V\cap T^*_ZX}\not=0$.
\end{lemma}
\begin{proof}
By the generic classification theorem in \cite{S-K-K},
$\she_Y^\infty\tens[\she_Y]\shn$ is a de Rham system
by shrinking $V$ if necessary.
There exists a smooth complex hypersurface $Z$ of $Y$ 
such that $T^*_ZY\subset V$.
Then by a quantized contact transform,
$\she_Y^\infty\tens[\she_Y]\shn$ and $H^1(\mu_Z(\sho_Y))$ are transformed to
$\bigl(\she_{\C^n}^\infty/(\sum_{i=1}^s\she_{\C^n}^\infty\partial_i)
\bigr)^{\oplus m}$ and
$H^1(\mu_{\{z_n=0\}}(\sho_{\C^n}))$ with $s<n$ and $m>0$.
In this case, the assertion is obvious.
\end{proof}

\medskip
\noindent
{\em Proof of \eqref{eq:b}}\quad
The proof is similar to that of \eqref{eq:a}.
Note that $V_\alpha$ contains $T^*_YX$.
By shrinking $Y$ if necessary, we may assume
that $\car(\shm)$ is equal to $V_\alpha$
on a neighborhood of a point $p$ in $T^*_YX$.
Then on a neighborhood of $p$,
$\she_X^\R\tens[\shd_X]\shm$
is isomorphic to
$\she_{X\hot Y}^\R\tens[\shd_Y]\shn$ for
a coherent $\shd_Y$-module $\shn$ with $\Supp(\shn)=Y$
on a neighborhood of $\pi(p)$.
One has then by \cite{S-K-K}
\eqn
H^{k}(\mu_Y(F))&\simeq&
\hom[\shd_X]\bl(\shm,H^k(\mu_Y(\sho_X))\br)\\
&\simeq&\hom[\she_X^\R]\bl(\she_X^\R\tens[\shd_X]\shm,H^k(\mu_Y(\sho_X))\br)\\
&\simeq&
\hom[\she_X^\R]\bl(\she_{X\hot Y}^\R
\tens[\shd_Y]\shn,H^k(\mu_Y(\sho_X))\br)\\
&\simeq&\pi^{-1}\hom[\shd_Y](\shn, \sho_Y)\vert_{T^*_YX}.
\eneqn
Hence the assertion follows from the following well-known result.
\end{proof}
\begin{lemma}
If a coherent $\shd_X$-module $\shm$ satisfies
$\Supp(\shm)=X$,
then one has $\Supp(\hom[\shd_X](\shm,\sho_X))=X$.
\end{lemma}

By the Riemann-Hilbert correspondence of
perverse sheaves and holonomic $\shd$-modules,
we have the following description of
the truncated microsupport
of perverse sheaves.

\begin{corollary}
Let $F\in D^b(\C_X)$ and
let $\{X_\alpha\}_{\alpha\in A}$ be a family of
complex submanifolds such that
$\ol{X_\alpha}$ and $\ol{X_\alpha}\setminus X_\alpha$
are closed complex analytic subsets
and $\Ss(F)=\bigcup_{\alpha\in A} \ol{T^*_{X_\alpha}X}$.
If $F$ is a perverse sheaf
{\rm(}i.e.\ there is a holonomic $\shd_X$-module $\shm$ such that
$F\simeq\rhom[\shd_X](\shm,\sho_X)${\rm)}, then one has
\eq\label{eq:perv}
\Ss_k(F)=\bigcup_{\codim X_\alpha\le k}
\ol{T^*_{X_\alpha}X}\quad\mbox{for any $k$.}
\eneq
Conversely if $F\in D^b(\C_X)$ is $\C$-constructible
and if it satisfies 
\eq\label{eq:pervs}
\Ss_k(F)\cup\Ss_k(\rhom(F,\C_X))\subset\smash{\bigcup_{\codim X_\alpha\le k}}
\ol{T^*_{X_\alpha}X}\quad\mbox{for any $k$,}
\eneq
then $F$ is a perverse sheaf.
\end{corollary}
\begin{proof}
It is a direct consequence of Theorem \ref{th:sskFcharM}
that the perversity of $F$ implies \eqref{eq:perv}.
Conversely assume \eqref{eq:pervs}.
In order to prove that $F$ is a perverse sheaf,
it is enough to show
that $F$ is microlocally isomorphic to
$\C_{X_\alpha}[-\codim X_\alpha]^{\oplus m}$
for some $m$ at a generic point of
$T^*_{X_\alpha}X$ by \cite[Theorem 10.3.12]{K-S1}.
By \cite{K-S1}, $F$ is isomorphic to $\C_{X_\alpha}[-\codim X_\alpha]\otimes K$
at a generic point of
$T^*_{X_\alpha}X$ for some $K\in D^b(\C)$.
Since $\mu_{X_\alpha}(F)$
must be in $D^{\ge \codim X_\alpha}(\C_{T^*_{X_\alpha}X})$
and $\mu_{X_\alpha}(F)\simeq
\C_{T^*_{X_\alpha}X}[-\codim X_\alpha]\otimes K$,
one has $K\in D^{\ge0}(\C)$.
Similarly,
$\mu_{X_\alpha}\bl(\rhom(F,\C_X)\br)\simeq
\C_{T^*_{X_\alpha}X}[-\codim X_\alpha]\otimes \rhom(K,\C)$
implies $K\in D^{\le0}(\C)$.
\end{proof}

{\small
Masaki Kashiwara\\
Research Institute for Mathematical Sciences,\\
Kyoto University, Kyoto 606-8502 Japan\\
masaki@kurims.kyoto-u.ac.jp
\vspace{5mm}

Teresa Monteiro Fernandes\\
Centro de {\'A}lgebra da Universidade de Lisboa, Complexo 2,\\
2 Avenida Prof. Gama Pinto, 1699 Lisboa codex Portugal\\
tmf@ptmat.lmc.fc.ul.pt
\vspace{5mm}

Pierre Schapira\\
Universit{\'e} Pierre et Marie Curie, case 82\\
Institut de Math{\'e}matiques\\
4, place Jussieu, 75252 Paris cedex 05 France\\
schapira@math.jussieu.fr\\
http://www.institut.math.jussieu.fr/{\~{}}schapira/
}


\begin{thebibliography}{15}

\bibitem {B-S}
J-M.~Bony and P.~Schapira,\,\,
{\em Existence et prolongement des solutions holomorphes
des \' equations aux d{\'e}riv{\'e}es partielles,}
\,Inventiones Math., \,\textbf{17},\, 95-105\, (1972).

\bibitem {E-K-S1}
P.~Ebenfelt, D.~Khavinson, H.~Shapiro, \, 
{\em Analytic continuation of Jacobi
polynomial expansions,} \, Indag. Math. S.N. \textbf{8(1)}\, 19-31 (1997).

\bibitem {E-K-S2}
\bysame, \, 
{\em Extending solutions of holomorphic
partial differential equations across real hypersurfaces}, 
\, J. London Math. Soc.
\textbf{(2) 57} 411-432 (1998).


\bibitem{Ho}
L.~H{\"o}rmander,
{\em The analysis of linear partial differential operators I,}
Grundlehren der Math. Wiss.\,\textbf{256}, Springer Verlag, (1983).


\bibitem {Ka-1}
M.~Kashiwara,
{\em Systems of Microdifferential Equations},
Progress in
Math.,\, \textbf{34},\, Birkh{\"a}user
(1983).

\bibitem{Ka-2}
\bysame,
{\em Algebraic Analysis,}\,
in Japanese, to be translated by the AMS.



\bibitem{K-S2}
M.~Kashiwara and P.~Schapira,
{\em Microlocal Study of Sheaves,}
\, Ast{\'e}risque \,\,\textbf{128},\, Soc. Math. France (1985).


\bibitem{K-S1}
\bysame,
{\em Sheaves on manifolds,}
Grundlehren der Math. Wiss. \textbf{292}, Springer Verlag (1990).

\bibitem{Le}
J.~Leray,
{\em Probl{\`e}me de Cauchy I,}\,
Soc. Math. France \,\textbf{85}, 389--430 (1957).


\bibitem{S}
P.~Schapira,
{\em Microdifferential systems in the complex domain},
Grundlehren der Math. Wiss., \textbf{269}, Springer Verlag (1985)

\bibitem{S-K-K}
M.~Sato, T.~Kawai and M.~Kashiwara,
{\em Hyperfunctions and pseudodifferential equations,}
Lecture Notes in Math., Springer \textbf{287}, \, 265--529 (1973).

\bibitem{To}
 F.~Tonin,
{\em Holomorphic extension for solutions of the characteristic Cauchy
Problem,}\,
Preprint (1998).

\bibitem{Ze}
M.~Zerner,
{\em Domaines d'holomorphie des fonctions v{\'e}rifiant une {\'e}quation
  aux d{\'e}riv{\'e}es partielles,} C.R. Acad. Sci.,\, \textbf{272},
\, 1646--1648 (1971).

\end{thebibliography}
\end{document}